\documentclass[12pt]{amsart}

\usepackage{amssymb}

\usepackage{anysize}
\marginsize{2cm}{2cm}{2cm}{2cm}

\renewcommand{\le}{\leqslant}

\usepackage{versions}
\usepackage{hyperref}
\hypersetup{
  colorlinks   = true, 
  urlcolor     = blue, 
  linkcolor    = blue, 
  citecolor   = red 
}
\usepackage{cancel}
\usepackage{yhmath}

\theoremstyle{definition}

\theoremstyle{remark}

\numberwithin{equation}{section}

{%
  \goodbreak
}
{%
  \goodbreak
}
{%
  ~]]\goodbreak\smallskip
}
\excludeversion{solution}
\excludeversion{answer}

\newcommand{\where}{\mathop{\ |\ }}

\renewcommand{\epsilon}{\varepsilon}
\renewcommand{\phi}{\varphi}
\renewcommand{\kappa}{\varkappa}
\renewcommand{\theta}{\vartheta}

\title{Huang's theorem and the exterior algebra}

\author{Roman~Karasev}

\thanks{Supported by the Federal professorship program grant 1.456.2016/1.4 and the Russian Foundation for Basic Research grants 18-01-00036 and 19-01-00169}

\address{Roman~Karasev, Moscow Institute of Physics and Technology, Institutskiy per. 9, Dolgoprudny, Russia 141700\newline \indent
Institute for Information Transmission Problems RAS, Bolshoy Karetny per. 19, Moscow, Russia 127994}
\email{r\_n\_karasev@mail.ru}
\urladdr{http://www.rkarasev.ru/en/}


\begin{document}

\maketitle

\setcounter{section}{1}

\subsection{Introduction/abstract}

In this note we give a version of Hao Huang's proof of the sensitivity conjecture, shedding some light on the origin of the magical matrix $A$ in that proof. For the history of the subject and the importance of this conjecture to the study of boolean functions, we refer to the original paper \cite{huang2019}. Here we only state the main result: Consider the boolean cube $Q_n=\{0,1\}^n$ as a graph, whose edges connect pairs of vertices differing in one coordinate. Then any its induced subgraph on greater than $2^{n-1}$ (the half) vertices has degree of some vertex at least $\sqrt{n}$.

\subsection{Exterior algebra and its linear endomorphism}

Let $V$ be a real vector space with basis $e_1,\ldots,e_n$, equip its dual $V^*$ with the dual basis $e_1^*,\ldots,e_n^*$. The exterior algebra $\wedge^* V$ consists of antisymmetric polylinear forms on $V$ and linear combinations of such forms of different degrees; its multiplication is the exterior product $\wedge$. Another useful operation is the interior product of $\omega\in \wedge^k V$ by a vector $v\in V$, given by
\[
(i_v \omega)(v_1,\ldots,v_{k-1}) = \omega(v, v_1, \ldots, v_{k-1}),
\]
this is a derivation of the exterior algebra. We refer to the textbook \cite{kostrman1989} for such algebraic basics.

Now choose a vector $v\in V$ and a linear form $\lambda\in V^*$ and consider the operator $A: \wedge^*V \to \wedge^*V$ defined as
\[
A(\omega) = i_v\omega + \lambda\wedge \omega.
\]
Let us determine its eigenvalues. Take the square and obtain
\[
A^2(\omega) = i_v(\lambda\wedge\omega) + \lambda\wedge i_v\omega = \lambda(v)\omega - \lambda\wedge i_v\omega + \lambda\wedge i_v\omega = \lambda(v) \omega.
\]
Hence $A^2$ is a scalar operator, and the eigenvalues of $A$ must be $\pm \sqrt{\lambda(v)}$. Since $A$ changes the parity of the degree, $A(\wedge^k(V)) \subseteq \wedge^{k-1}(V)\oplus\wedge^{k+1}(V)$, its trace is zero and therefore its eigenvalues $\sqrt{\lambda(v)}$ and $-\sqrt{\lambda(v)}$ have the same multiplicity. Since its square is nonzero scalar, $A$ must also be semi-simple and there is a splitting
\[
\wedge^*(V) = G_A^+ \oplus G_A^-
\]
into the positive and negative eigenspaces of $A$ of dimension $2^{n-1}$ each, provided $\lambda(v) > 0$ (otherwise we would need to complexify).

\subsection{Basis of the exterior algebra and the boolean cube}

We may choose the standard basis in $\wedge^*(V)$, consisting of products $e^*_{i_1} \wedge \dots \wedge e^*_{i_k}$,
for sequences $1\le i_1 < \dots < i_k \le n$. Such sequences may be indexed by $0/1$-vectors of length $n$, having $1$ precisely at positions $i_1,\ldots,i_k$. Hence we may assume that the basis of the exterior algebra corresponds to the vertices of the boolean cube $Q_n=\{0,1\}^n$, we may take $\beta\in Q_n$ and consider its corresponding basis element $e^*_\beta\in \wedge^*(V)$.

The operation $i_v$, when applied to a basis form, may only change one $1$ to $0$ in the boolean notation and linearly combine the results of such changes. The operation $\lambda\wedge$ may only do the opposite, linearly combine changes of one $0$ to $1$. Hence, for the linear operation $A$ from the previous section, $A(e^*_\beta)$ expresses as a linear combination of $e^*_\gamma$, where $\gamma$ are adjacent to $\beta$ in the $1$-skeleton graph of the cube. Let us denote the adjacency relation by $\beta\leftrightarrow\gamma$.

Passing to the sensitivity conjecture \cite{huang2019} about subgraphs of $Q_n$, we choose more than $2^{n-1}$ vertices $H\subseteq Q_n$ and consider the corresponding linear subspace $H\subset \wedge^*(V)$, denoted by the same letter. From dimension considerations $H$ must have nonzero intersection with $G_A^+$, take some nonzero $\omega \in G_A^+\cap H$. Then we have the relation:
\[
A(\omega) = \sqrt{\lambda(v)}\omega.
\]
Writing in the basis $\omega = \sum_{\beta\in H} \omega_\beta e^*_\beta$ and choosing the largest $|\omega_\beta|$ (let us flip $\omega$ to make its coordinate $\omega_\beta$ positive), we then have in coordinates
\[
\sqrt{\lambda(v)} \omega_\beta = \sum_{\gamma\leftrightarrow\beta} A_{\beta\gamma} \omega_\gamma,
\]
where $A_{\beta\gamma}$ are the matrix elements of $A$ in the chosen basis. If $\beta$ contains one more $1$ compared to $\gamma$ (denote this by $\gamma\rightarrow\beta$), then $|A_{\beta\gamma}| = |\lambda_k|$, where $k$ is the position of this extra $1$ and $\lambda_k$ is the corresponding coordinate of $\lambda$. In the opposite case, $\beta\rightarrow\gamma$, $|A_{\beta\gamma}| = |v_\ell|$, where $\ell$ is the position of the change and $v_\ell$ is the coordinate of $v$ from the definition of $A$. Therefore
\[
\sqrt{\lambda(v)} \omega_\beta = \left|\sum_{\gamma\leftrightarrow\beta} A_{\beta\gamma} \omega_\gamma\right| 
\le \|\lambda\|_\infty  \sum_{\gamma\rightarrow\beta} |\omega_\gamma| + \|v\|_\infty \sum_{\beta\rightarrow\gamma} |\omega_\gamma| 
\le  \|\lambda\|_\infty \sum_{\gamma\rightarrow\beta,\ \gamma\in H} \omega_\beta + \|v\|_\infty \sum_{\beta\rightarrow\gamma,\ \gamma\in H} \omega_\beta,
\]
where $\|\cdot\|_\infty$ denotes the maximal coordinate of a vector or a linear form. Dividing by $\omega_\beta$, we obtain
\[
\sqrt{\lambda(v)} \le \|\lambda\| \#\{\gamma\in H \where \gamma\rightarrow\beta\} + \|v\|_\infty \#\{\gamma\in H\where \beta\rightarrow\gamma\}.
\]

Since we are interested in lower bounds on the degree of a vertex of the graph induced by $H$, we need to increase $\lambda(v)$ given $\|\lambda\|_\infty$ and $\|v\|_\infty$. It then makes sense to put $\lambda_k \equiv a$ and $v_\ell \equiv b$, for $a,b>0$, thus obtaining
\[
\sqrt{nab} \le a \#\{\gamma\in H\where \gamma\rightarrow\beta\} + b \#\{\gamma\in H\where \beta\rightarrow\gamma\}.
\]
The case $a=b=1$ delivers the sensitivity conjecture, the degree of a vertex in the induced by $H$ subgraph is bounded from below by $\sqrt{n}$. The proof in \cite{huang2019} is essentially an elementary translation of this argument.

Other values of $a$ and $b$ may provide some information, when we view $Q_n$ as an oriented graph and want to estimate its incoming and outgoing degrees of a vertex. Dividing by $\sqrt{ab}$ and substituting $C = \sqrt{\frac{a}{b}}$, we obtain
\[
\sqrt{n} \le C \#\{\gamma\in H\where \gamma\rightarrow\beta\} + \frac{1}{C} \#\{\gamma\in H\where \beta\rightarrow\gamma\}
\]
for any $C>0$. Of course, the vertex of $H$ where this inequality holds depends on $C$.

\subsection*{Acknowledgment}
The author thanks Arseniy Akopyan for useful remarks on this exposition.

\bibliography{../Bib/karasev}
\bibliographystyle{abbrv}
\end{document}